\documentclass [12pt]{article}
\usepackage{graphicx,amssymb,amsfonts,latexsym,amsmath,amsthm,times}
\usepackage{epsfig}
\usepackage{fancyhdr}
\usepackage{color}
\usepackage[english]{babel}
\newtheorem{theorem}{Theorem}

\newtheorem{corollary}[theorem]{Corollary}

\newtheorem{lemma}[theorem]{Lemma}

\title{ Comparison of Perron and Floquet eigenvalues 
 in age structured cell division cycle models}
 \date{}
\author{Jean Clairambault$^{a,b}$,
  St\'ephane Gaubert$^{c,d}$   and Thomas Lepoutre$^{a,e}$ \footnote{Corresponding
author. E-mail: thomas.lepoutre@inria.fr}}
\begin{document}
\maketitle

\centerline{$^a$ INRIA, projet BANG, Domaine de Voluceau, BP 105, 78156 Le Chesnay Cedex France}

\centerline{$^b$ INSERM U 776,
H\^opital Paul-Brousse,
14, Av. Paul-Vaillant-Couturier
F94807 Villejuif cedex}

\centerline{$^c$ INRIA Saclay -- \^Ile-de-France, projet MAXPLUS} 
\centerline{$^d$ CMAP, Ecole Polytechnique, 91128 Palaiseau Cedex, France}
\centerline{$^e$ UPMC Univ Paris 06, UMR 7598, Laboratoire Jacques-Louis Lions, F-75005, Paris, France}

\vspace*{1cm}
\renewcommand{\thefootnote}{}
\footnotetext{This work was partially supported by the European Union,
through the Network of Excellence BIOSIM (contract number LSHBCT-2004-005137)
and the Scientific Targeted Project TEMPO (contract number LSHG-ct-2006-037543),
and by the RFBR-CNRS joint grant number 05-01-02807}
\renewcommand{\thefootnote}{\arabic{footnote}}
\abstract{
We study the growth rate of a cell population  that follows an 
age-structured PDE with time-periodic coefficients. Our motivation 
comes from the comparison between experimental tumor growth curves 
in mice endowed with intact or disrupted circadian clocks, 
known to exert their influence on the cell division cycle.
We compare the growth rate of the model controlled by a time-periodic 
control on its coefficients with the growth rate of stationary 
models of the same nature, but with averaged coefficients. 
We firstly derive a delay differential equation which allows us to prove 
several inequalities and equalities on the growth rates. We also 
discuss about the necessity to take into account  the structure of 
the cell division cycle for  chronotherapy modeling. Numerical simulations illustrate the results.
}

\textbf{Key words:} cell cycle, circadian rhythms, chronotherapy, structured PDEs, delay differential equations.

{\bf AMS subject classification:} 35F05, 35P05, 35P15, 92B05, 92D25.


\vspace*{1cm}

\setcounter{equation}{0}
\section{Cell cycle control and circadian rhythms}
The cell division cycle is the process by which the eukaryotic cell duplicates its DNA content and then divides itself in two daughter cells. 
This process is normally controlled by various physiological 
mechanisms that ensure homeostasis of healthy tissues, that control genome integrity (e.g.\ cyclins and cdks, p53, repair enzymes, etc.), launching programmed cell death (apoptosis)
if the DNA is irreversibly damaged (see \cite{OMorgan} for a complete presentation). The system of control has been extensively studied and modeled (see e.g.\  
\cite{Goldbetermitotic,Keener,Novak}
or \cite{Tyson}) using ordinary differential equations.
The cell division can be modeled through branching processes (see \cite{ArinoKimmel}), integral equations, delay differential equations
(see \cite{SB}) and also
many structured PDE models (for an overview, see \cite{Arino,ArinoSanchez,MetzDiekmann})
 where the structuring variables can be age (\cite{Murray1}), size (\cite{Perthamebook}) or more recently cyclin content 
(\cite{Bekkal1,Bekkal2,Doumic}).

Most living organisms exhibit circadian rhythms (from Latin \textit{circa diem}, ``roughly a day'') 
which allow them to adapt to an environment that varies with a periodicity of 
24h. These rhythms  can be observed even in the smallest  biological functional unit, the cell. 
The problem we are studying is the growth of cell populations (undergoing the cell division cycle described above) under the pressure of circadian rhythms.
Circadian rhythm effects on the cell cycle turn out to be important in tumor proliferation. This is observed by several experiments involving
a major disruption of circadian rhythms in mice. In these experiments 
it can be seen that the growth of tumors is significantly enhanced
 in mice in which the pacemaker circadian clock has been drastically 
 perturbed, either through neurosurgery, or through light-dark cycle 
 disruption (see e.g.\ \cite{Filipski2002,Filipski2005}). Moreover, in 
 the clinic, taking advantage of the influence exerted by  circadian 
 clocks on anticancer drug metabolism and on the cell division 
 cycle has led in the past 15 years to successful applications in the {\it chronotherapy of 
 cancers}, particularly colorectal cancer (see \cite{Levi}). This motivates modeling the circadian 
rhythm in simple cell cycle models and studying these effects on the growth rate of a cell population.
\begin{figure}[htbp]
\centerline{\includegraphics[scale=1]{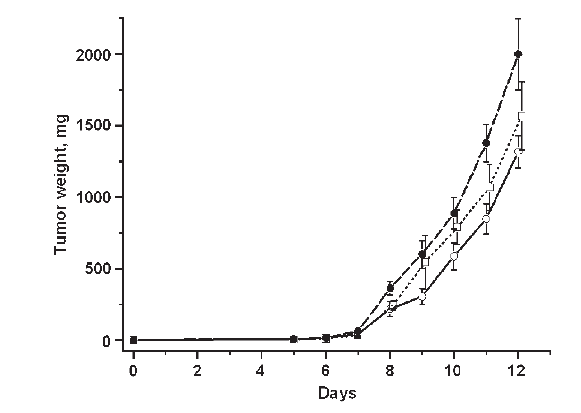}}
\caption{Effects of the perturbation of light-dark cycle on tumor proliferation (reproduced from \cite{Filipski2005}). In  clock-perturbed mice 
(black dots), the tumor 
proliferates much faster than in control mice (white dots).(By courtesy of Elizabeth Filipski).}
 \label{Filipski}
\end{figure}

Contrary to our first idea, the growth rate of a cell population 
described by a physiologically structured PDE model with time-periodic 
control is not necessarily lower than in a model of the same nature, but with 
a time-averaged control \cite{CGP,CMP,CMP2}. 

The goal here is twofold. Firstly we analyze how modeling assumptions lead to define various growth rates 
under the effects of circadian rhythms. Secondly we  model the effect of chronotherapy on these growth rates.

In the second section we recall the definition of these various growth rates, in terms of Perron and Floquet eigenvalues of a linear Von Foerster-
Mc-Kendrick model. We also discuss known inequalities between them.
In the third section we study a simple division 
model, for which we establish (in Theorem~\ref{Comparaison-locale}) strict inequalities
comparing the growth rate in the stationary (Perron) and periodic (Floquet)
cases. These inequalities are proved by studying a related time delay system
(which is similar to the one considered in~\cite{SB}). This model is used to confirm the impossibility 
to derive a general comparison between the Perron and Floquet eigenvalues defined in the second section. In the fourth section, we give an argument for using multiphase models to represent chronotherapy, taking better into account the structure of the cell cycle and particularly the existence of various phases.
We provide numerical simulations to illustrate our results.
In a first appendix, we give the detailed proof of the existence of the
solution of the eigenproblem, by applying the Krein-Rutman theorem. In a second
appendix, we derive analytical formul\ae\ for the eigenelements in a specific
multiphase case, which yield further information on their behavior and can be
used to validate numerical experiments.
\setcounter{equation}{0}
\section{The  model}
%
%
%
%
%
%

\subsection{The renewal equation}

 We base our study on a cell population that follows the classical renewal equation structured in age with periodic coefficients representing
the effect of circadian rhythms 
\begin{equation}
\label{renewal}
\left\lbrace\begin{array}{l}
\frac{\partial  }{\partial t}n(t,x)+\frac{\partial  }{\partial x}n(t,x)+d(t,x)n(t,x)=0,\\[0.3cm]
n(t,x=0)=\int_0^\infty B(t,x)n(t,x)dx.
\end{array}\right.
\end{equation}
Here $n(t,x)$ represents the density of cells of age $x$ in the cycle at time $t$, $d(t,x), B(t,x)$ represent respectively the 
death rate,
and the birth rate. Both  these coefficients are $T$-periodic in time. 
We define the growth rate of the population in terms of an eigenproblem. The growth rate $\lambda_F$ 
(F for Floquet as for ODEs with periodic coefficients) is defined as the unique real number 
$\lambda_F$, such that there is a solution $N$ to  the problem
\begin{equation}
\label{Floquet-renewal}
\left\lbrace\begin{array}{l}
\frac{\partial  }{\partial t}N(t,x)+\frac{\partial  }{\partial x}N(t,x)+[\lambda_F+d(t,x)]N(t,x)=0,\\[0.3cm]
N(t,x=0)=\int_0^\infty B(t,x)N(t,x)dx,\\[0.3cm]
N>0,\qquad T-\text{periodic}.
\end{array}\right.
\end{equation}
We refer to \cite{MMP} for conditions of existence for $\lambda_F$ (and to the appendix for the case of division models).

\subsection{Comparison of eigenvalues}

We use the following notations. For a $T$-periodic function $f$ we define,
\begin{eqnarray*}
\langle f\rangle &=&\frac{1}{T}\int_0^T f(t)dt \qquad \text{the arithmetical average,}\\
\langle f\rangle_g &=& \exp\Big(\frac{1}{T}\int_0^T \log f(t)dt\Big) \qquad \text{the geometrical average, when } f>0.
\end{eqnarray*}
%
%
%
%
%
%

It may seem natural to introduce the following stationary problem (Perron eigenproblem), in which the death and birth rates are averaged
\begin{equation}
\left\lbrace\begin{array}{l}
\label{Perron-renewal}
\frac{d}{dx}N_P(x) +\lbrack \lambda_P+\langle d(x)\rangle\rbrack N_P(x)=0, \\[0.3cm] 
N_P(0)=\int_0^\infty \langle B(x)\rangle N_P(x)dx=1,\\[0.3cm]
N_P(x)>0.
\end{array}\right.
\end{equation}

It is shown in \cite{CMP,CMP2} that, when $B$ does not depend on time, the inequality $\lambda_F\geq\lambda_P$ holds. In the present paper, we show that this inequality does not carry over to the case of a time dependent $B$. It should be noted, however, that there is a general inequality, established in \cite{CGP}, which relates $\lambda_F$ with the solution of the following eigenproblem in which an arithmetical average of the death rate is taken, whereas the geometrical average of the birth rate is taken,
\begin{equation}
\left\lbrace\begin{array}{l}
\label{Geom-renewal}
\frac{d}{dx}N_g(x) +\lbrack \lambda_g+\langle d(x)\rangle\rbrack N_g(x)=0, \\[0.3cm] 
N_g(0)=\int_0^\infty \langle B(x)\rangle_g N_g(x)dx=1,\\[0.3cm]
N_g(x)>0.
\end{array}\right.
\end{equation}

\begin{theorem}[{\cite{CGP}}]\label{ComparaisonFloquetPerronGeo}
The eigenvalues defined in (\ref{Floquet-renewal}) and (\ref{Geom-renewal}) satisfy
\begin{equation*}
\lambda_F\geq\lambda_g.
\end{equation*}
\end{theorem}
This result suggests that there is no general inequality between $\lambda_P$ and $\lambda_F$,  because the inequality which follows from convexity is $\lambda_F\geq \lambda_g$. Moreover, it follows from the standard arithmetico-geometrical inequality, 
\begin{equation*}
\lambda_P\geq\lambda_g.
\end{equation*} 
  Such a general comparison cannot hold between $\lambda_F$ and $\lambda_P$, 
as shown in the next section. To go further we use a more specific model.
\setcounter{equation}{0}
\section{A simple one-phase division model}
%
%
%
%
%
%

\subsection{Model and main results}

We model the cell cycle with the following PDE which is a particular case of (\ref{renewal}), 
\begin{equation*}
\label{1Malthus}
\left\lbrace\begin{array}{l}
\frac{\partial}{\partial t}n(t,x)+\frac{\partial}{\partial x}n(t,x) +[d(t)+ K_0\psi(t) \chi_{[a,+\infty[}(x) ]n(t,x)=0, \\[0.3cm]
n(t,0)=2K_0\psi(t) \int_{a}^\infty n(t,x)dx,
\end{array}\right.
\end{equation*}
where $K_0>0$ is a constant, $\psi>0$ is a $T$-periodic function with 
\begin{equation}
\langle \psi\rangle=1.\label{moypsi}
\end{equation} The term
 $K_0\psi(t)\chi_{[a,+\infty[}$ represents the division rate, $d(t)$ is the apoptosis rate (we assume it to be $T$-periodic). We have denoted by 
$\chi_E$ the indicator function of set $E$. Finally, $\psi(t)$ represents a nonnegative periodic control exerted on division.
As before we look for the growth rate $\lambda_F$ of such a system. It is defined so  that there is a solution  to the Floquet eigenproblem,
\begin{equation}\label{Floquet1}
\left\lbrace\begin{array}{l}
\frac{\partial}{\partial t}N(t,x)+\frac{\partial}{\partial x}N(t,x) +\big[ \lambda_F+d(t)+ K_0\psi(t)\chi_{[a,+\infty[}(x) \big]N(t,x)=0, \\[0.3cm]
N(t,0)=2K_0\psi(t) \int_{a}^\infty N(t,x)dx,\\[0.3cm]
N>0, \; \text{$T$-periodic},
\end{array}\right.
\end{equation}
and we normalize $N$ by
\begin{equation*}
\int_0^T\int_0^\infty N(t,x)dxdt=1.
\end{equation*}

As we already know a general comparison result for the geometrical eigenvalue $\lambda_g$ defined in (\ref{Geom-renewal}), we are now only interested in the 
comparison of $\lambda_F$ and $\lambda_P$, the latter quantity defined by requiring the existence of a solution to the Perron eigenproblem already defined in (\ref{Perron-renewal}) which
here reads
\begin{equation}
\left\lbrace\begin{array}{l}
\label{ARIT1}
\frac{d}{d x}N_P(x) +\lbrack \lambda_P+\langle d\rangle+K_0\chi_{[a,+\infty[}(x)\rbrack N_P(x)=0, \\[0.3cm] 
N_P(0)=2K_0\int_{a}^\infty N_P(x)dx,\\[0.3cm]
N_P>0,
\end{array}\right.
\end{equation} 
and we normalize $N_P$ by
\begin{equation*}
N_P(0)=2K_0\int_{a}^\infty N_P(x)dx=1.
\end{equation*}

We are interested in evaluating the effect of the periodic control $\psi(t)$ on the growth of the system.
Therefore we denote by $\lambda_F(a,\psi)$ and by $\lambda_P(a)$ the above defined eigenelements so as to keep track of 
the problem parameters.\\

The following theorem implies that there is no possible general comparison between $\lambda_F$ and $\lambda_P$.
\begin{theorem}\label{Comparaison-locale}
For all continuous positive $T$-periodic functions $\psi$ satisfying (\ref{moypsi}), we have
\begin{equation}\lambda_F(a=T,\psi)=\lambda_P(T)=\lambda_F(a=T,1),\label{equality}\end{equation}
and for $a$ in  a neighborhood of T, we have, provided $\psi\not\equiv 1$
\begin{equation*}
\begin{array}{c}
\lambda_F(a,\psi)>\lambda_P(a)=\lambda_F(a,1)\quad \text{for} \quad a<T,\\[0.3cm]
\lambda_F(a,\psi)<\lambda_P(a)=\lambda_F(a,1)\quad \text{for} \quad a>T.
\end{array}\label{inequalities}
\end{equation*}
\end{theorem}

The proof of this theorem is presented in the next sections. The computations done in section 4.1 insure that, without loss of generality, we can suppose $d\equiv0$.

\begin{figure}[htbp]
\centerline{\includegraphics[scale=0.5]{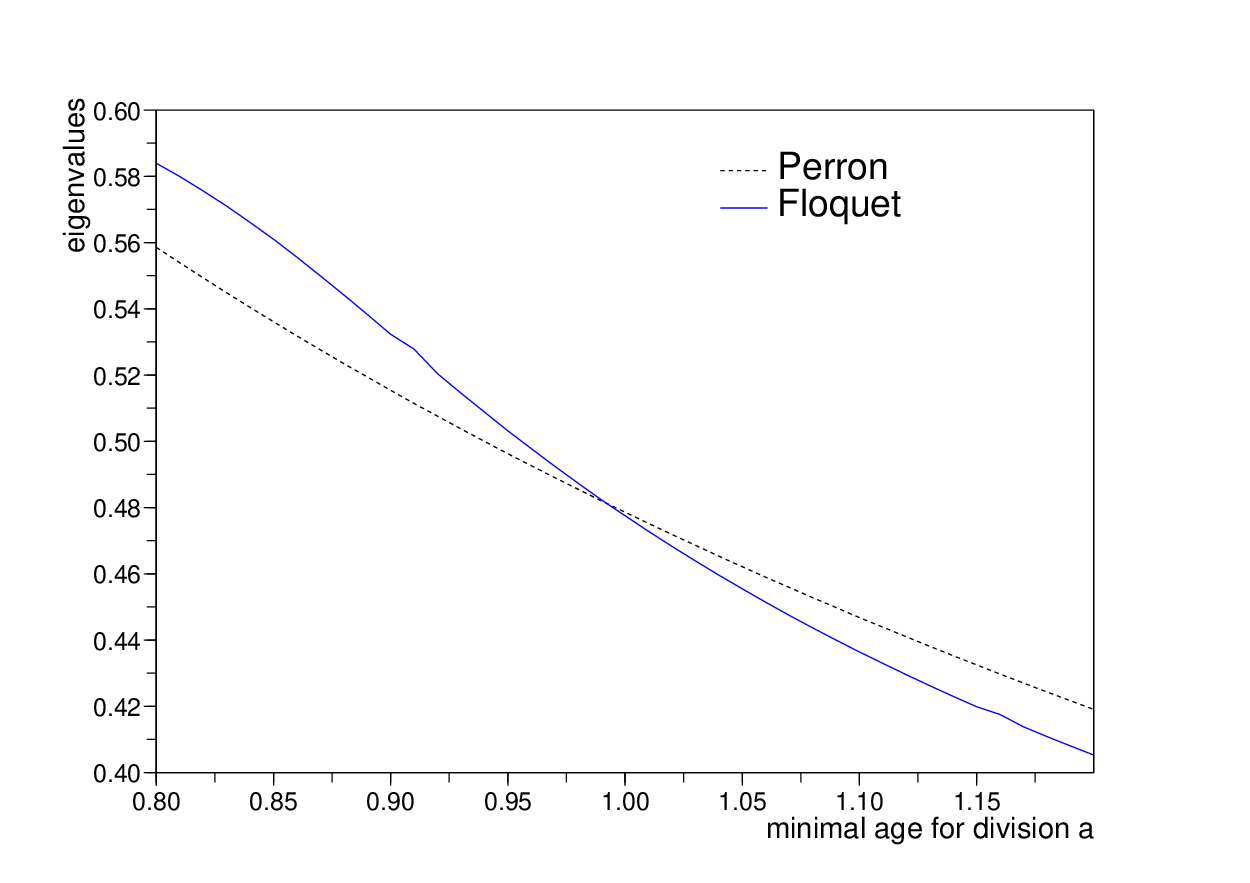}}
\caption{Crossing of the Perron and Floquet curves (detail) for $\psi=\psi_{\sin}$.}
 \label{fig:Floquet-Perron_detail}
\end{figure}

Numerical results are presented in figures \ref{fig:Floquet-Perron_detail} and \ref{fig:Floquet-Perron_global} which illustrate this theorem. Graphically, for fixed $\psi$, this predicts firstly that the curves of $\lambda_F(a,\psi)$ (Floquet curve) and $\lambda_P(a)$ (Perron curve) must cross each other for $a=T$, secondly that the Floquet curve should be above the Perron curve before (i.e., for $a<T$) the crossing and below this curve after it (i.e., for $a>T$). A possible interpretation is that for a better adaptation 
(in the sense of higher proliferation), the cell cycle should be shorter than $24$h; an effect already observed  in \cite{SB}.

\subsection{Proof of Theorem 2, part 1 (a delay differential equation)}

Throughout the proof, we use the shorter notations $\lambda_F$ and $\lambda_P$ instead of $\lambda_F(a,\psi)$ and 
$\lambda_P(a)$ when there is no possible confusion.\\
To find more information on $\lambda_F$ we derive a delay differential equation.\\
We integrate (\ref{Floquet1}) with respect to age over $[a,\infty[$.  We get
\begin{eqnarray*}
\frac{d}{dt}\int_{a}^\infty N(t,x)dx +N(t,\infty)-N(t,a)+[\lambda_F
+K_0\psi(t)]\int_{a}^\infty N(t,x)dx=0.
\end{eqnarray*}

{From} the formula of characteristics and the boundary condition in (\ref{Floquet1}), 
\begin{equation*}
\begin{array}{l}
N(t,a)=N(t-a,0)e^{-\lambda_Fa},\\
N(t,a)=2K_0e^{-\lambda_F a}\psi(t-a)\int_{a}^\infty N(t-a,x)dx.
\end{array}
\end{equation*}

We set $P(t)=\int_{a}^\infty N(t,x)dx $. Since we have $N(t,\infty)=0$ (see the appendix) we obtain the delay differential equation
\begin{eqnarray}
\label{EDOret1}
\dot P(t)+\Big( \lambda_F+K_0\psi(t)\Big) P(t)
=2K_0\psi(t-a)P(t-a)e^{- \lambda_F a}.
\end{eqnarray}

\subsection{Proof of Theorem 2, part 2 (equality of growth rates for $a=T$)}

The comparison between $\lambda_P$ and $\lambda_F$ is based on the following formula for $\lambda_P$.


\begin{lemma}\label{SOLPERRON}The Perron eigenvalue defined in (\ref{ARIT1}) satisfies
\begin{equation}\forall a > 0,\qquad \frac{\lambda_P+K_0}{2K_0}e^{\lambda_P a}=1.\label{solPerron}\end{equation}
\end{lemma}


\noindent {\bf Proof.} {From} (\ref{ARIT1}), we have, for $x\geq a$, $N_P(x)=e^{-(\lambda_P+K_0)x+K_0a}.$
We insert that in the boundary condition and obtain
$$1=2K_0\int_{a}^\infty e^{-(\lambda_P+K_0)x+K_0a}dx,$$
$$1=2K_0\frac{1}{\lambda_P+K_0}e^{- \lambda_Pa}.$$\qed

\begin{corollary}\label{lambdaPpositif}
The Perron eigenvalue defined in (\ref{ARIT1}) satisfies
\begin{equation*}
\forall a>0, \qquad \lambda_P>0.
\end{equation*}
\end{corollary}

\noindent {\bf Proof.} This follows from Lemma \ref{SOLPERRON} and the remark 
$$\forall a>0,\forall \lambda\leq 0,\qquad \frac{\lambda+K_0}{2K_0}e^{\lambda a}\leq \frac{1}{2}.$$\qed

To obtain (\ref{equality}), we divide (\ref{EDOret1}) by $P$ and find 
\begin{eqnarray*}
\frac{\dot P(t)}{P(t)} =
- \lambda_F-K_0\psi(t)+2K_0\psi(t-a)\frac{P(t-a)}{P(t)}e^{-\lambda_F a}.
\end{eqnarray*}
When we take the average over a period, we get (since $P$ is $T$-periodic in time by its definition as $N$ is)
$$0=-(\lambda_F+K_0) +2K_0e^{-\lambda_F a} \bigg\langle \psi(t-a)\frac{P(t-a)}{P(t)}\bigg\rangle ,$$
\begin{equation}\label{fonclambda}
\frac{\lambda_F+K_0}{2K_0}e^{\lambda_F a }= \bigg\langle \psi(t-a)\frac{P(t-a)}{P(t)} \bigg\rangle.
\end{equation}
Now we consider the particular case $a=T$. As $P$ is $T$-periodic $P(t-a)=P(t)$. Hence, for $a=T$, we arrive at
\begin{equation}
\frac{\lambda_F+K_0}{2K_0}e^{\lambda_F a }= \bigg\langle \psi(t-a)\frac{P(t-a)}{P(t)}\bigg\rangle
=\langle \psi\rangle
=1.\label{equalitybis}
\end{equation}
This equality is the same for $\lambda_F$ as the one described in lemma 1 for $\lambda_P$. As we know that the mapping

\begin{eqnarray*}
\lambda &\mapsto & \frac{\lambda+K_0}{2K_0}e^{\lambda a}.
\end{eqnarray*}
is increasing  on $[-K_0  ,+\infty[$ from $0$ to $+\infty$ and is negative elsewhere, 
there is only one solution to (\ref{equalitybis}) which is also given by (\ref{solPerron}) and the result (\ref{equality}) is proved. \qed

\subsection{Proof of Theorem 2, part 3 (local comparison around $a=T$)}

We fix $\psi \not\equiv 1$. We study the variations of $\frac{\lambda_F+K_0}{2K_0}e^{\lambda_F a}$
around $a=T$. 
{From} (\ref{fonclambda}), we know:
$$\frac{\lambda_F+K_0}{2K_0}e^{\lambda_Fa}=
\bigg\langle \psi(t-a)\frac{P(t-a)}{P(t)}\bigg\rangle=\bigg\langle \psi(t)\frac{P(t)}{P(t+a)}\bigg\rangle,$$
therefore
\begin{eqnarray*}
\frac{\partial}{\partial a} \frac{\lambda_F+K_0}{2K_0}e^{\lambda_Fa} &=& \frac{\partial}{\partial a}
\bigg\langle \psi(t)\frac{P(t)}{P(t+a)}\bigg\rangle, \\
&=& \bigg\langle \psi(t) \frac{\partial P}{\partial a}(t)\frac{1}{P(t+a)}\bigg\rangle  
+ \bigg\langle \psi(t)\frac{-P(t)}{P^2(t+a)}\frac{\partial }{\partial a}\big(P(t+a)\big) \bigg\rangle. 
\end{eqnarray*}
Recalling that $P$ depends on $a$ (as $N$ and $\lambda_F$ do), we have 
\begin{equation*}
\frac{\partial  }{\partial a} P(t+a) =\frac{\partial P}{\partial a} (t+a)+\dot P(t+a).
\end{equation*}

We then split the computations
\begin{eqnarray*}
\frac{\partial}{\partial a} \frac{\lambda_F+K_0}{2K_0}e^{\lambda_Fa}&=&
\bigg\langle \psi(t) \frac{\partial P}{\partial a}(t)\frac{1}{P(t+a)}\bigg\rangle  
-\bigg\langle \psi(t)\frac{P(t)}{P^2(t+a)}\bigg(\frac{\partial P}{\partial a} (t+a)+\dot P(t+a)\bigg)\bigg\rangle,\\
&=&\bigg\langle \psi(t)\frac{1}{P(t+a)}\bigg(\frac{\partial P}{\partial a}(t)-
\frac{P(t)}{P(t+a)}\frac{\partial P}{\partial a}(t+a)\bigg)\bigg\rangle\\
&&-\bigg\langle \psi(t)\frac{P(t)}{P^2(t+a)}\dot P(t+a)\bigg\rangle.
\end{eqnarray*}

For $a=T$, the first term vanishes, and $P(t+a)=P(t)$ i.e., 
\begin{equation*}
\frac{\partial}{\partial a} \frac{\lambda_F+K_0}{2K_0}e^{\lambda_Fa}=-\bigg\langle \psi(t)\frac{P(t)}{P^2(t)}\dot P(t)\bigg\rangle
=-\bigg\langle \psi(t)\frac{\dot P(t)}{P(t)}\bigg\rangle.
\end{equation*}
To compute this we  again make use of the ODE (\ref{EDOret1}) which we multiply by $\dfrac{\psi}{P}$
\begin{eqnarray*}
\psi(t)\frac{\dot P(t)}{P(t)}=- \lambda_F\psi(t)-K_0\psi^2(t)
+2K_0\psi(t-a)\psi(t)\frac{P(t-a)}{P(t)}e^{- \lambda_F a}.
\end{eqnarray*}
Averaging on a period we still get, for $a=T$, 
\begin{eqnarray*}
\bigg\langle \psi(t)\frac{\dot P(t)}{P(t)}\bigg\rangle=-\lambda_F-K_0\langle\psi^2\rangle
+2K_0\langle\psi^2\rangle e^{- \lambda_F a}.
\end{eqnarray*}
Using (\ref{equalitybis}),  we arrive at
\begin{eqnarray*}
\bigg\langle\psi(t)\frac{\dot P(t)}{P(t)}\bigg\rangle&=& 
-\lambda_F  -K_0\langle\psi^2\rangle+
\langle\psi^2\rangle(\lambda_F+K_0)=\lambda_F(\langle \psi^2\rangle -1).
\end{eqnarray*}
We now have  the derivative at $a=T$,
\begin{equation}
\label{deriv}
\frac{\partial}{\partial a}_{\mid_{a=T}} \frac{\lambda_F+K_0}{2K_0}e^{\lambda_F a} =
-\lambda_F(\langle \psi^2\rangle -1).
\end{equation}

We use here the  notations $\lambda_F'(T)$ for $\frac{\partial \lambda_F}{\partial a}_{\mid_{a=T}} $ and $\lambda_F(T)=\lambda_P(T)$
to recall that we are studying the local behavior of $\lambda_F$ and $\lambda_P$ around $a=T$, ($\psi$ is fixed).  We  can directly compute
\begin{equation*}
\frac{\partial}{\partial a}_{\mid_{a=T}} \frac{\lambda_F+K_0}{2K_0}e^{\lambda_F a}=\lambda_F'(T)\frac{e^{\lambda_F(T) T}}{2K_0}
+(\lambda_F'(T)T+\lambda_F(T))\frac{K_0+\lambda_F(T)}{2K_0}e^{\lambda_F(T) T},
\end{equation*}
Therefore, using (\ref{equality}) and (\ref{solPerron}), we obtain
\begin{equation*}
\frac{\partial}{\partial a}_{\mid_{a=T}} \frac{\lambda_F+K_0}{2K_0}e^{\lambda_F a} = \lambda_F'(T)
\bigg(\frac{e^{\lambda_P(T) T}}{2K_0}+T\bigg)+\lambda_F(T), 
\end{equation*}
so that, using (\ref{equality}) and (\ref{deriv}), we have  
\begin{equation*}
\lambda_F'(T)=
\frac{-\lambda_P(T)\langle\psi^2\rangle}{T+\frac{e^{\lambda_P(T) T}}{2K_0}}. 
\end{equation*}
Similarly we have
\begin{equation*}
\lambda_P'(T)=
\frac{-\lambda_P(T)}{T+\frac{e^{\lambda_P(T) T}}{2K_0}}. 
\end{equation*}
Therefore,
\begin{equation*}
\lambda_P'(T)-\lambda_F'(T)=\frac{\lambda_P(T)(\langle\psi^2\rangle -1)}{T+\frac{e^{\lambda_P(T) T}}{2K_0}}. 
\end{equation*}
Thanks to corollary \ref{lambdaPpositif}, $\lambda_P(T)$ is positive. The assumption (\ref{moypsi}) leads to
$$\langle\psi^2\rangle-1=\big\langle (\psi -1)^2\big\rangle >0.$$
Finally we obtain 
\begin{equation}
\lambda_P'(T)-\lambda_F'(T)>0,\label{diffderiv}
\end{equation}
and the second statement of the  theorem follows then immediately from  (\ref{equality}) and (\ref{diffderiv}).\qed 

%
%
%
%
%
%
\setcounter{equation}{0}
\section{Modeling chronotherapy}
In the following we propose a model for chronotherapy by the introduction of a periodic death rate due to the effect of a drug on our cell division cycle model.
\subsection{Limit of single-phase division models}\label{sec:Limit of single-phase division models}

We consider a population of cells following a general division equation with apoptosis rate $d$. As above, all  coefficients are $T$-periodic with respect to time.
\begin{equation*}
\left\lbrace\begin{array}{l}
\frac{\partial}{\partial t}n(t,x)+\frac{\partial}{\partial x}n(t,x) +\big(d(t,x)+ K(t,x) \big)n(t,x)=0, \\[0.3cm]
n(t,0)=2\int_0^\infty K(t,x)n(t,x)dx.
\end{array}\right.
\end{equation*}
We consider the Floquet eigenproblem associated with this equation
\begin{equation*}
\left\lbrace\begin{array}{l}
\frac{\partial}{\partial t}N(t,x)+\frac{\partial}{\partial x}N(t,x) +\big(d(t,x)+ K(t,x)+\lambda_F \big)N(t,x)=0, \\[0.3cm]
N(t,0)=2\int_0^\infty K(t,x)N(t,x)dx,\\[0.3cm]
N>0,\quad \int_0^T\int_0^\infty N(t,x)dxdt=1
\end{array}\right.
\end{equation*}
We propose to model the effect of chronotherapy by adding a time $T$-periodic, age-independent death rate $\gamma(t)$ representing the effect of a drug (for instance we may consider 
$\gamma$ proportional to the quantity of drug in the body). The cell population now follows  the equation
\begin{equation*}
\left\lbrace\begin{array}{l}
\frac{\partial}{\partial t}n(t,x)+\frac{\partial}{\partial x}n(t,x) +[d(t,x)+ K(t,x) +\gamma(t)]n(t,x)=0, \\[0.3cm]
n(t,0)=2\int_0^\infty K(t,x)n(t,x)dx.
\end{array}\right.
\end{equation*}
The Floquet eigenproblem for this equation reads 
\begin{equation*}
\left\lbrace\begin{array}{l}
\frac{\partial}{\partial t}N^{\gamma}(t,x)+\frac{\partial}{\partial x}N^\gamma(t,x) +\big(d(t,x)+ K(t,x)+\gamma(t)+\lambda_F^\gamma \big)N^\gamma(t,x)=0, \\[0.3cm]
N^\gamma(t,0)=2\int_0^\infty K(t,x)N^\gamma(t,x)dx,\\[0.3cm]
N^\gamma>0,\quad T-\text{periodic} \quad \int_0^T\int_0^\infty N^\gamma(t,x)dxdt=1.
\end{array}\right.
\end{equation*}
\begin{lemma}\label{lemma:lambdad}
The Floquet eigenvalue $\lambda_F^\gamma$ defined above satisfies
\begin{equation*}
\lambda_F^\gamma=\lambda_F-\langle \gamma\rangle.\end{equation*}
\end{lemma}

\noindent {\bf Proof.} We define $\tilde \gamma=\gamma-\langle \gamma\rangle$, $ \Gamma(t)=\int_0^t \tilde \gamma(s)ds$. Noticing that $\Gamma$ is $T$-periodic,
we define the function $M$ by $M(t,x)~=~N(t,x)e^{\tilde \Gamma(t)}$. It satisfies
\begin{equation*}
\left\lbrace\begin{array}{l}
\frac{\partial}{\partial t}M(t,x)+\frac{\partial}{\partial x}M(t,x) +\big(d(t,x)+ K(t,x)+\gamma(t)+\lambda_F-\langle\gamma\rangle \big)M(t,x)=0, \\[0.3cm]
M(t,0)=2\int_0^\infty K(t,x)M(t,x)dx,\\[0.3cm]
M>0,\quad T-\text{periodic}.
\end{array}\right.
\end{equation*}
Therefore $\lambda_F^\gamma=\lambda_F-\langle\gamma\rangle$ and up to a renormalization $M=N^\gamma$.\qed

This result expresses that with such a simple model, chronotherapy is inefficient, since changing the moment of administration of a drug (in symbols, changing $\gamma(t)$ into $\gamma(t+\theta)$ where $\theta$ is a real number) has no effect on the growth rate. In other words, in such one-phase models, this effects depends on $\langle \gamma\rangle$. Only the total daily dose of the drug is relevant!

\subsection{Using multiphase models}

We now consider more realistic multiphase models. We use the additional ingredient that the real cell division cycle is multiphasic because of the existence of checkpoints between phases (mainly at the  G1/S and G2/M transitions) at which it can be arrested if genome integrity is not preserved. We consider a cell cycle model with $I$ phases where $I>1$ (for instance $I=4$ if we want to represent the classical phases  G1-S-G2-M). We study $I$ populations of cells, $n_i(t,x)$ being the density of cells of age $x$ in phase $i$ at time $t$. We use the convention $I+1=1$
\begin{equation}\label{multi}
\left\lbrace\begin{array}{l}
\frac{\partial}{\partial t}n_i(t,x)+\frac{\partial}{\partial x}n_i(t,x)+[K_{i\rightarrow i+1}(t,x)+d_i(t,x)]n_i(t,x)=0,\\[0.3cm]
n_{i+1}(t,0)=\int_0^\infty K_{i\rightarrow i+1}(t,y) n_i(t,y) dy, \qquad 1<i\\[0.3cm]
n_{1}(t,0)=2\int_0^\infty K_{I\rightarrow 1}(t,y) n_I(t,y) dy,\\[0.3cm]
n_i(0,x)=n^0_i(x) \quad given.
\end{array}\right.
\end{equation}
Here $K_{i\rightarrow i+1}$ represents the transition rate from phase $i$ to $i+1$. At the end of phase $I$ division occurs with rate $K_{I\rightarrow 1}$. To be as general as possible, we have considered death rates $d_i$ in phase $i$. As above, the coefficients are time $T$-periodic and we can consider the Floquet eigenproblem
\begin{equation}\label{multiFloquet}
\left\lbrace\begin{array}{l}
\frac{\partial}{\partial t}N_i(t,x)+\frac{\partial}{\partial x}N_i(t,x)+[K_{i\rightarrow i+1}(t,x)+d_i(t,x)+\lambda]N_i(t,x)=0,\\[0.3cm]
N_{i+1}(t,0)=\int_0^\infty K_{i\rightarrow i+1}(t,y) N_i(t,y) dy, \qquad 1<i\\[0.3cm]
N_{1}(t,0)=2\int_0^\infty K_{I\rightarrow 1}(t,y) N_I(t,y) dy,\\[0.3cm]
N_i>0,\quad T-\text{periodic},\quad \sum\limits_i \int_0^1\int_0^\infty N_i dxdt=1.
\end{array}\right.
\end{equation}
We also consider  the adjoint eigenproblem
\begin{equation}\label{multiFloquetadjoint}
\left\lbrace\begin{array}{l}
\frac{\partial}{\partial t}\phi_i(t,x)+\frac{\partial}{\partial x}\phi_i(t,x)-[K_{i\rightarrow i+1}(t,x)+d_i(t,x)+\lambda]\phi_i(t,x)=K_{i\rightarrow i+1}\phi_{i+1}(t,0),\\[0.3cm]
\frac{\partial}{\partial t}\phi_I(t,x)+\frac{\partial}{\partial x}\phi_I(t,x)-[K_{I\rightarrow 1}(t,x)+d_I(t,x)+\lambda]\phi_I(t,x)=2K_{I\rightarrow 1}\phi_{1}(t,0),\\[0.3cm]
\phi_i>0,\quad T-\text{periodic},\quad \sum\limits_i \int_0^\infty N_i\phi_idxdt=1.
\end{array}\right.
\end{equation}
To model the effect of chronotherapy, we consider a cytotoxic drug acting only on a specific phase (for instance 5-Fluorouracil acts on S-phase, see \cite{LACG} for instance and the references therein) and, as in the previous section we represent its action by an additional death rate in phase $j$, $\gamma(t)$ (we replace in phase $j$ $d_j$ by $d_j+\gamma$) . We also define eigenelements for the modified equation $(\lambda^\gamma,N^\gamma,\phi^\gamma)$. We multiply the first line  of (\ref{multiFloquet}) (version with $d_j$ replaced by $d_j+\gamma$, $N_i$ by $N_i^\gamma$ and $\lambda$ by $\lambda^\gamma$) by $\phi_i$, and (\ref{multiFloquetadjoint}) by $N_i^\gamma$. Summing over $i$ and integrating over age and time, 
we obtain 
\begin{equation}\label{Ngamma}
(\lambda-\lambda^\gamma)\sum\limits_i \int_0^1\int_0^\infty N^\gamma_i \phi_i dxdt=\int_0^1\gamma(t)\int_0^\infty N^\gamma_j \phi_j dxdt.
\end{equation}
We shall not have here the problem encountered with one-phase models. We study the effect of a death rate $\gamma(t+\theta)$. We denote $\lambda^{\varepsilon,\theta},N^{\varepsilon,\theta}$ the eigenelements associated to an additional death rate $\varepsilon\gamma(t+\theta)$ in phase $j$. We define $F(\varepsilon,\theta)$ by 
\begin{equation}\label{F}
F(\varepsilon,\theta)=\lambda-\lambda^{\varepsilon,\theta}=\frac{\int_0^1\varepsilon\gamma(t+\theta)\int_0^\infty N^{\varepsilon,\theta}_j \phi_j dxdt}{\sum\limits_i \int_0^1\int_0^\infty N^{\varepsilon,\theta}_i \phi_i dxdt}.
\end{equation}
As we have $\lambda=\lambda^{0,\theta}$ for any $\theta$, $F(0,\theta)\equiv 0$. Particularly it does not depend on $\theta$. The question is: does $F(\varepsilon,\theta)$ depend on $\theta$ for fixed $\varepsilon$? To assess this question, we compute using  dominated convergence
\begin{equation}
\frac{\partial \lambda(\varepsilon,\theta)}{\partial \varepsilon}\mid_{\varepsilon=0}=\lim\limits_{\varepsilon\rightarrow 0} \dfrac{F(\varepsilon,\theta)}{\varepsilon}=\int_0^1\gamma(t+\theta)\int_0^\infty N_j\phi_jdxdt.
\end{equation}
Therefore if neither the function $\gamma(.)$ nor the function $\int_0^\infty N_j\phi_j(.,x)dx$ are constant (contrarily to one-phase models, there are no compensating effect making $\int_0^\infty N_j\phi_j(.,x)dx$ constant, see for instance the computations of the appendix), then $\lim\limits_{\varepsilon\rightarrow 0} \dfrac{F(\varepsilon,\theta)}{\varepsilon}$  depends on $\theta$ (we mean it is not a constant function of $\theta$) and so is (at least for small $\varepsilon$) $F(\varepsilon,.)$. In this case the Taylor first order approximation around $0$ of $\lambda$: $\lambda(\varepsilon,\theta)\approx\lambda+\varepsilon\int_0^1\gamma(t+\theta)\int_0^\infty N_j\phi_jdxdt$ is not a constant function of $\theta$ and neither is $\lambda(\varepsilon,\theta)$, at least for small values of $\varepsilon$. We illustrate this property numerically in the next section (see figure \ref{fig:chronotherapie1}). It seems that the Taylor first order approximation is a very good approximation of the growth rate for a reasonable range of values of the amplitude $\varepsilon$.

\setcounter{equation}{0}
\section{Numerical simulations}\label{sec:Numerical results}
%
%
%
%
%
%
We illustrate the theorems proved above by several numerical simulations. We firstly present the numerical scheme, then we give several algorithmic properties. 
Finally tests are presented.
\subsection{Discretization}
In our numerical simulations we consider a pure division model  :
\begin{equation}
\label{Nonapop1}
\left\lbrace\begin{array}{l}
\frac{\partial}{\partial t}n(t,x)+\frac{\partial}{\partial x}n(t,x)+K_0\psi(t)\chi_{[a,+\infty[}(x)n(t,x)=0, \\[0.3cm]
n(t,0)=2K_0\psi(t)\int_{a}^\infty n(t,x)dx.
\end{array}\right.
\end{equation}

Consider time and age increments $\Delta t,\Delta x$ and denote by $\kappa_i$ and $\psi^k$,
the quantities $\kappa_i =K_0\chi_{[a,+\infty[}(i\Delta x)$ and $\psi^k=\psi(k\Delta t)$. 
Choosing first order finite differences, we obtain from equation (\ref{Nonapop1}) 
the following approximation with an error of order $O(\left|\Delta t\right| + \left|\Delta x\right|)$ 

\[ \frac{n_i^{k+1}-n_i^k}{\Delta t}+\frac{n_i^k-n_{i-1}^k}{\Delta x}+
\kappa_{i}\psi^{k+1} n_{i}^{k+1}=0,\quad 1\leq i\leq I,\]
where $\{0\dots I\}$ is the set of all values of $i$ to be considered in the discretization.
Taking $\Delta t=\Delta x$ ($CFL=1$), we  obtain the following compact discretization scheme:
\begin{equation}
\label{disc}
\left\lbrace\begin{array}{l}
n_{i}^{k+1} = \frac{n_{i-1}^{k}}{1+\Delta t \kappa_i \psi^{k+1}},\quad 1\leq i \leq I, \\[0.3cm]
n_{0}^{k+1} = 2\psi^{k}\sum\limits_{0\leq i\leq I} \kappa_i n_{i}^{k}\Delta t.
\end{array}\right.
\end{equation}

Assume $\psi$ is periodic of period $T\geq 0$ and consider a grid over $[0,T]\times[0,I\Delta t]$, 
consisting of squares with sides of length $\Delta t= T/N_{T}$, for some $N_T\in\mathbb{N}$ (and $I$  large
 enough, particularly $I\Delta t>a$ and $I+1>N_T$). 
Then, the populations at time $(k+1)\Delta t$ for all ages in $[0,I\Delta t]$ can be obtained from the corresponding populations 
at time $k \Delta t$ as follows:
\begin{equation}\label{disc3}
\left(\begin{array}{c} 
n_{0}^{k+1}\\ 
n_{1}^{k+1}\\ 
\vdots \\ 
n_{I}^{k+1} 
\end{array}\right) = \left(\begin{array}{cccc}
\scriptstyle{2\psi^{k}\kappa_0\Delta t} & \ldots & \scriptstyle{2\psi^{k}\kappa_{_{I - 1}}\Delta t} & \scriptstyle{2\psi^{k}\kappa_{_{I}}\Delta t} \\
\frac{1}{1+\Delta t \psi^{k+1} \kappa_1} & \ldots & 0 & 0 \\
\vdots & \ddots & \vdots & \vdots \\
0 & \ldots & \frac{1}{1+\Delta t \psi^{k+1}\kappa_{_{I}}} & 0
\end{array}\right)\left(\begin{array}{c} 
n_{0}^{k}\\ 
n_{1}^{k}\\ 
\vdots \\ 
n_{I}^{k} 
\end{array}\right)
\end{equation}
 
It is clear that the matrix in (\ref{disc3}) depends only on the time index  $k$ and is periodic of period $N_T$. 
We denote $M_k$ this matrix and the vectors respectively $n^k$ and $n^{k+1}$. The equation (\ref{disc3}) can be written $n^{k+1}=M_kn^k$.

\subsection{Approximating the eigenvalue}

The algorithm has already been discussed in \cite{Emilio1phase}.
We recall that the growth rate is defined as the unique real $\lambda_F$ such that (\ref{Nonapop1}) admits solutions of the form $N(t,x)e^{\lambda_F t}$
with $N>0$ and $N(.,x)$ is periodic. We can approximate it thanks to:
\begin{lemma}[Discrete Floquet theorem]\hfill

There exists a unique real $\lambda$ and a unique sequence of vectors $ (\mathcal{N}^k)_{k\in\mathbb{N}}$ ,
$\mathcal{N}^k=\Big(\mathcal{N}^k_i\Big)_{0\leq i\leq I}$
 such that
\begin{equation} 
\mathcal{N}_i^k >0,\qquad \sum\limits_{i=0}^{I}\mathcal{N}_i^0=1,\label{pos}
\end{equation}
\begin{equation}
k\mapsto (\mathcal{N}^k)  \quad \text{is} \quad N_T\text{-periodic},\label{per}
\end{equation}
\begin{equation}
n^k, \quad \text{defined by }n^k=\mathcal{N}^k e^{\lambda .k \Delta t} \quad  \text{is solution to } 
(\ref{disc3})\label{solution}
\end{equation}
\end{lemma}
\noindent {\bf Proof.}
The proof is standard and we recall it for the sake of completeness. It is 
based on the Perron Frobenius theorem. 
First we prove uniqueness. Supposing there exists such $n^k$, we have
\begin{eqnarray}
n^1&=&M_0n^0 \nonumber,\\
n^2&=&M_1n^1=M_1M_0n^0 \nonumber,\\ 
&\ldots& \\ 
n^{k+1}&=&M_kn^k =M_kM_{k-1}\dots M_1M_0n^0 ,\\ 
&\ldots& \nonumber\\
 n^{N_T}&=& M_{N_T-1}M_{N_T-2}\cdots M_1 M_{0}n^{0}.\label{MatriceM}
\end{eqnarray}
We define
$$\mathbb{M}=M_{N_T-1}M_{N_T-2}\cdots M_1 M_{0},$$
thus, (\ref{MatriceM}) reads $n^{N_T}=\mathbb{M}n^0$.
\qed
\begin{lemma} The matrix
$\mathbb{M}$ is nonnegative
 and primitive (and therefore is irreducible). 
\end{lemma}
\noindent {\bf Proof.}
The nonnegativity is obvious. To prove the primitivity, the key point is $I+1>N_T$ and $I\Delta t\geq a+2\Delta t$.
For some $\varepsilon>0$ we have for any $k$, if we denote by $\operatorname{Id}_k$ the identity matrix of order $k$,
$$M_k\geq \varepsilon \left(\begin{array}{cr}
0\ldots 0  &1\;  1\\
\operatorname{Id}_{I} &0
\end{array}\right)=\varepsilon W.$$
Notice that $W$ is the Wielandt matrix of order $I+1$ which is known to be primitive (see \cite{HornJohnson}). 
Therefore for some $p$, $W^p>0$ and thus for  
$q N_T\geq p$,
$$\mathbb{M}^q\geq \varepsilon^{q N_T}W^{q N_T}>0,$$ which yields the primitivity of $\mathbb{M}$, the spectral radius of which, denoted here by $\rho$ is then positive.
We denote by $\rho$
its spectral radius. We have $\rho>0$.\\

Back to the proof of the discrete Floquet theorem, we have 
 $$n^{N_T}=e^{\lambda N_T\Delta t}\mathcal{N}^{N_T}  =e^{\lambda T}\mathcal{N}^{0}=e^{\lambda T} n^0.$$
Hence we have $\mathbb{M}n^0=e^{\lambda T} n^0$. This means that $n^0$ is a positive eigenvector of $\mathbb{M}$ 
associated to a positive eigenvalue $e^{\lambda T}$.
{From} the Perron-Frobenius theorem, $e^{\lambda T} =\rho$ and 
$n^0=\mathcal{N}^0$ is the (unique) associated eigenvector. The solution is unique.\\
Conversely, if we know the Perron eigenvector $V$ and the Perron eigenvalue $\rho$ of $\mathbb{M}$, then
the sequence $\Big(\mathcal{N}^k\Big)^{k\in \mathbb{N}}$ defined by 

$$\left\lbrace\begin{array}{l}
\mathcal{N}^0=V,\\[0.3cm]
\mathcal{N}^{k+1}=e^{-\lambda \Delta t}.M_k \mathcal{N}^{k},
\end{array}\right.$$

 satisfies (\ref{pos}),(\ref{per}) and (\ref{solution}) for $\lambda=\log(\rho(\mathbb{M}))$. \qed
 
For multiphase models, the idea is mainly the same.  
To compute $\rho=e^{\lambda T}$ the spectral radius of $\mathbb{M}$, 
the \textit{power algorithm} is used. It converges thanks to the primitivity of $\mathbb{M}$.
\subsection{Numerical results}
\begin{figure}[htbp!]
\centerline{\includegraphics[scale=0.5]{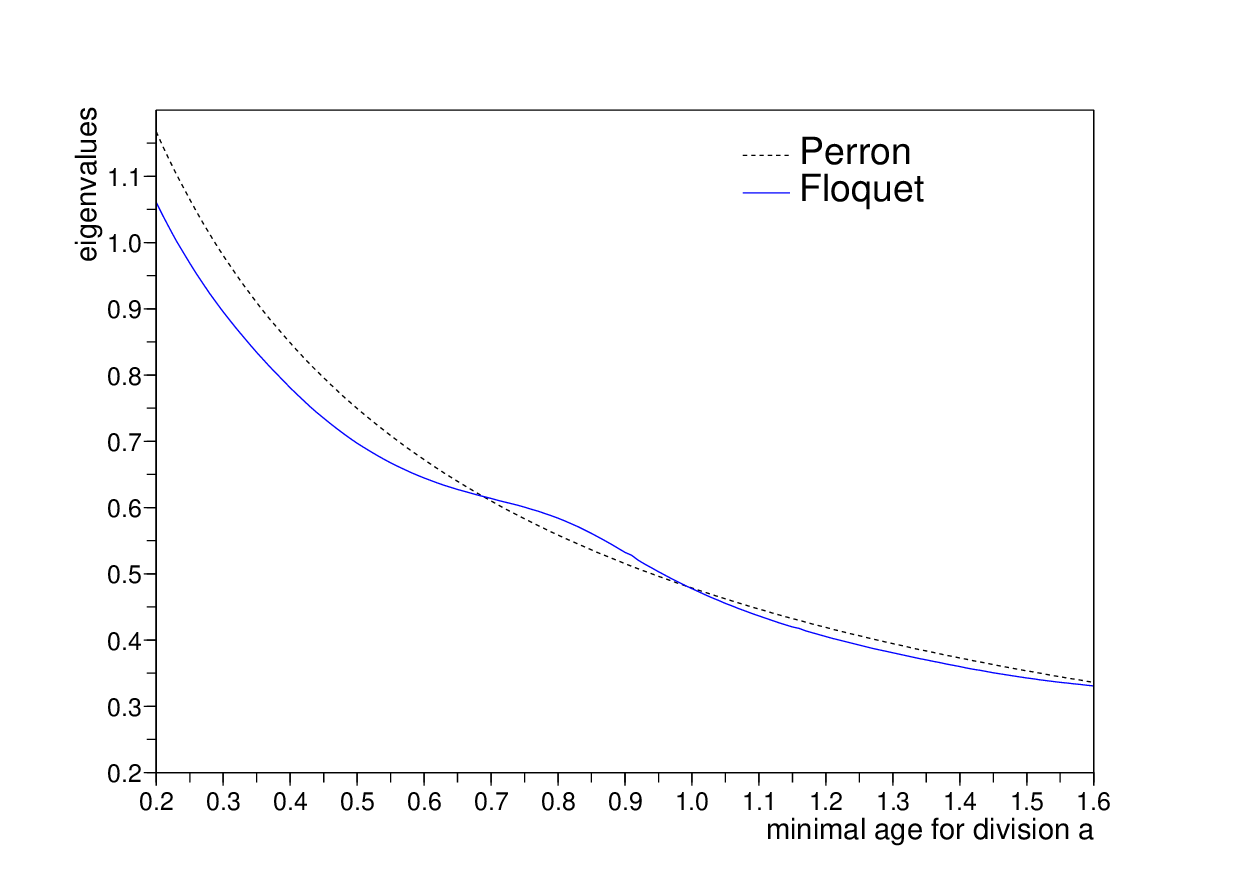}}
\caption{Crossing of the Perron and Floquet curves  for $\psi=\psi_{\sin}$.}
 \label{fig:Floquet-Perron_global}
\end{figure}

First we present some numerical results to illustrate theorem \ref{Comparaison-locale}. We scale $T=1$. We fix the value of $K_0$ to $2$ and test various 
periodic function $\psi$. We plot the curves
\begin{eqnarray*}
a&\rightarrow &\lambda_F(a,\psi),\\
a&\rightarrow &\lambda_P(a)
\end{eqnarray*}
We recall that the eigenvalues for the Perron problem can be directly computed thanks to lemma \ref{SOLPERRON}. {From} theorem \ref{Comparaison-locale}, we know that these curves cross for $x$-coordinate $a=T$, the second part of the theorem tells us that we expect (locally) the curve for $\lambda_F$ to be above the curve for $\lambda_P$ for $a<T$ and below it for $a>T$. 
We plot the curves $\lambda=\lambda_P(a)$ and $\lambda=\lambda_F(a,\psi)$ for our functions $\psi$ and look 
at the crossing of curves around $T$ (on the simulations, $T=1$). We  also give a more global view of $\lambda_F(a,\psi_{sin})$ and $\lambda_P$ in figure \ref{fig:Floquet-Perron_global} 
to illustrate the fact that the comparison is only local. 
\begin{table}[!h]
\begin{center}
\begin{tabular}{|c|c|c|}
\hline
Name of the function & Formulation on the interval $[0,1[$ & $\langle \psi^2\rangle$ \\\hline 
$\psi_{sq}$ (square wave) & $1.8\chi_{[0,1/2[}(t)+0.1\chi_{[1/2,1[}$ & 1.81 \\[0.1cm]\hline
&&\\[0.1mm]
$\psi_{pk}$ (peak function) & $0.1+ht/\delta \chi_{[0,\delta[}(t)+(2h-ht/\delta )\chi_{[\delta,2\delta[}(t)$ & 1.99 \\[0.1cm]\hline
$\psi_{sin}$ (sinusoidal) & $1+0.9\cos(2\pi t)$ & 1.405 \\\hline
\end{tabular}
\caption[periodic functions]{Functions $\psi$ for the simulations}
\label{periodic_functions}
\end{center}
\end{table}
Here, the parameters $h$ and $\delta$ are respectively set to $3$ and $0.3$.

\begin{figure}[htbp!]
\centerline{\includegraphics[scale=0.5]{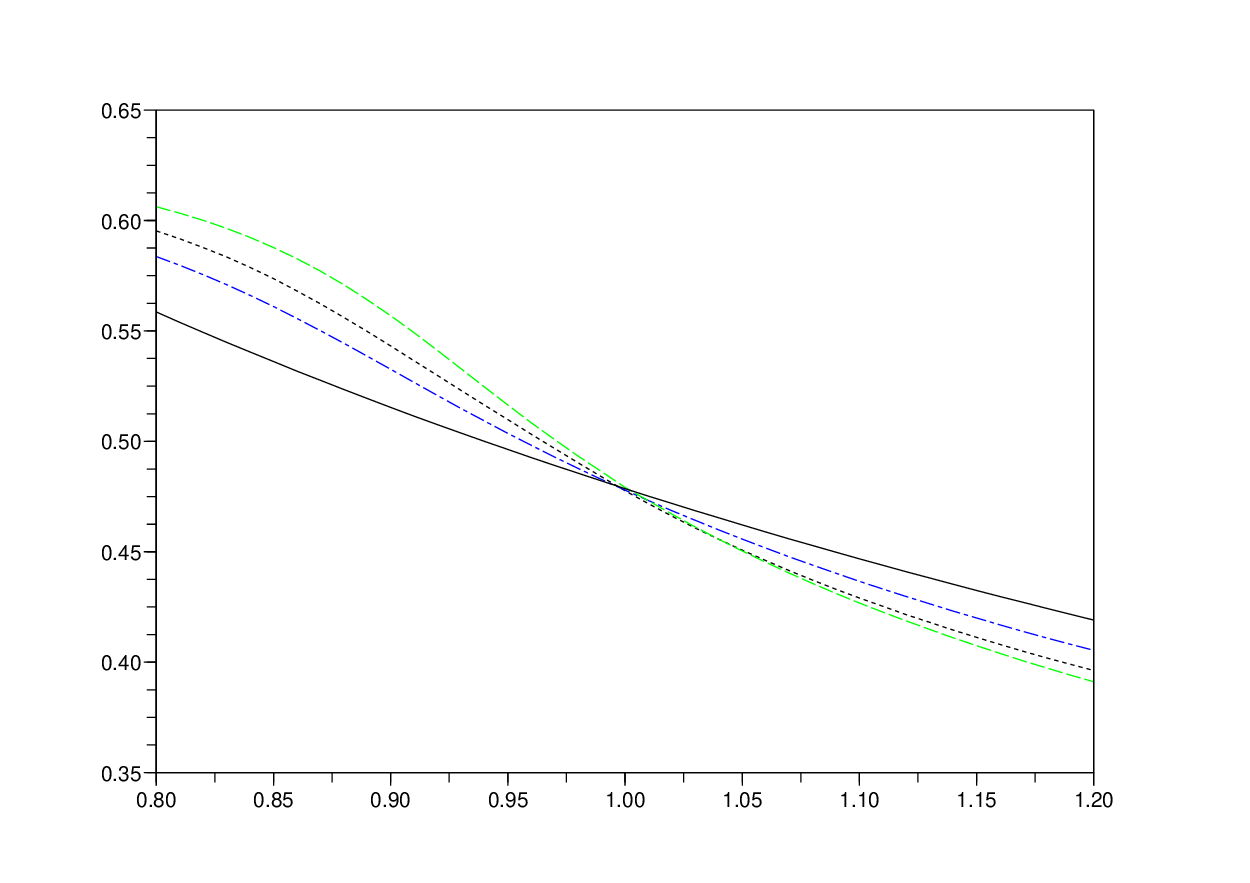}}
\caption{Crossing of the Perron and Floquet curves  for $\psi=\psi_{\sin}$ (dash dot), $\psi_{sq}$ (dots) and $\psi_{pk}$ (long dash).}
 \label{fig:Floquet-Perron_multi}
\end{figure}

{From} the last part of the demonstration of theorem \ref{Comparaison-locale}, we expect,
$$\frac{\partial\lambda_F(a,\psi_{pk})}{\partial a}_{\mid_{a=T}}>\frac{\partial\lambda_F(a,\psi_{sq})}{\partial a}_{\mid_{a=T}}>
\frac{\partial\lambda_F(a,\psi_{sin})}{\partial a}_{\mid_{a=T}},$$
we give figure \ref{fig:Floquet-Perron_multi} as a confirmation. Finally we give some simulations to illustrate our remarks on chronotherapy.

For the chronotherapy simulation we use the following parameter: we fix $I=3$ (we consider S and G2 as a single phase). The parameter $\gamma$ is a periodic function (with strong variations on a period to have a stronger effect of the parameter $\theta\in(0,1)$). We compute the eigenvalue for a death rate in phase 2 (phase S-G2) having the value $\varepsilon \gamma(t+\theta)$. We test several value of $\varepsilon$ to determine whether or not the amplitude of the death rate changes the relative behavior of the eigenvalue with respect to $\theta$.The coefficients have the form:
$$K_{i\rightarrow i+1}(t,x)=K_i\psi_i(t)\chi_{[a_i,\infty[}(x),$$
where $K_i,a_i$ are  positive, $\psi_i$ is a positive $1$ periodic function.
We give a simulation for the case described in the appendix (a case for which we can compute explicitly $\int_0^\infty N_2\phi_2(t,x)dx$). We fix $K_i=10$ for all $i$, $a_1=10/24, a_2=12/24=0.5, a_3=2/24$, $\psi(t)=1+0.9\cos(2\pi t)$ and $\psi_i$ defined from $\psi$ as in the appendix. We choose $\gamma(t)=\cos^6(2\pi t)$. With these choices of coefficients, we compute
$$\int_0^\infty N_2\phi_2(t,x)dx=C-C'\sin(2\pi t),$$
where $C$ and $C'$ are positive constants. Therefore,
$$\lim_{\varepsilon\rightarrow 0} \frac{\lambda^{\varepsilon,\theta}-\lambda^0}{\varepsilon}=C+C'\sin(2\pi\theta).$$
\begin{figure}[htbp!]
\centerline{\includegraphics[scale=0.5]{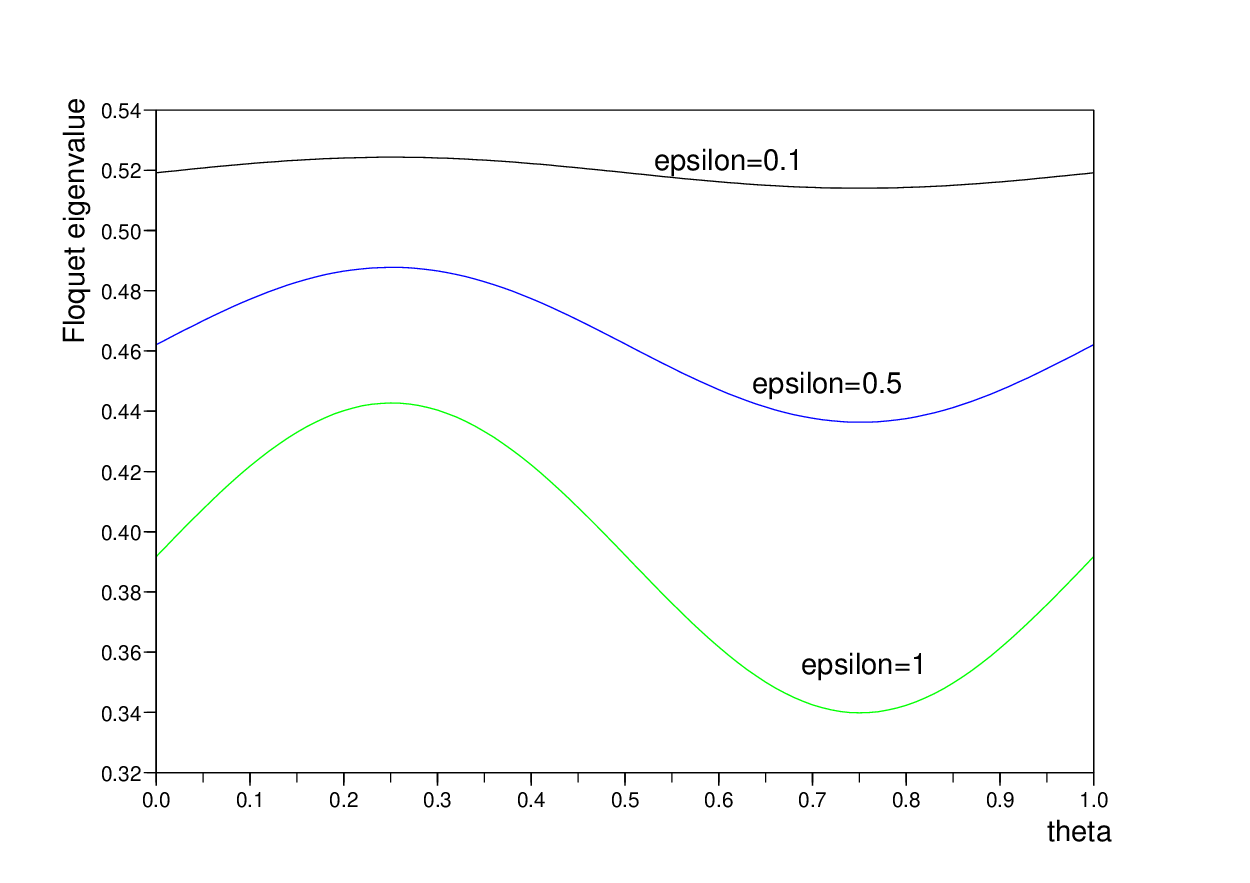}}
\caption{Variation of the Floquet eigenvalue with respect to the parameter $\theta$ for various amplitude for fixed $\gamma$ and amplitude $\varepsilon=0.1,0.5,1$ (from left to right).}
 \label{fig:chronotherapie1}
\end{figure}

In figure \ref{fig:chronotherapie1}, we remark especially that the location of the optimal phase does not depend on $\varepsilon$ (since we have $\theta_{\text{optimal}}=\frac{1}{4}$ whatever the value of $\varepsilon$) and corresponds exactly to the value of $\theta$ maximizing $\sin(2\pi\theta)$, i.e., minimizing $\int_0^1\gamma(t+\theta)\int_0^\infty N_2\phi_2dxdt$.
\newpage
 
%
%
%
%
%
%

\section*{Concluding remarks}
The results of the present paper show that the periodic control on the transition rate $K_{i\to i+1}$ of
cell cycle models yields richer behaviors than in the case in which only the death rates $d_i$ are subject
to a periodic control~\cite{CMP,CMP2}. In particular, the inequality of~\cite{CMP,CMP2} does not carry over. This
is, to our knowledge, the first time that such results are shown -on special cases of the control- analytically,
thus confirming numerical results first shown in~\cite{CMP,CMP2}.

Our results also indicate that multiphase cell proliferation models are the simplest candidates to represent
the effects of chronotherapy. Indeed, as shown in section~\ref{sec:Limit of single-phase division models}, in
single-phase models, in the simple case when only death rates $d_i$ are controlled by a periodic
forcing term, the growth rate $\lambda$ is modified by a term depending only on the average over a period
of the forcing term, so that no phase of the periodic control function can be relevant to account for differences
in the resulting growth rate, contrary to what is observed in chronotherapy~\cite{Levi}. Furthermore such multiphase models take into account the existence of multiple checkpoints , and we know from cell cycle physiology that the minimal number of checkpoints to consider is 2: at G1/S and G2/M.

We performed numerical and graphical results of section~\ref{sec:Numerical results}, on a 3-phase model with
1-periodic control on all phase transition functions $K_{i\to i+1}$, where one represents
chronotherapy as a 1-periodic death term $\varepsilon\gamma$ of amplitude $\varepsilon$ acting on the
second phase (S/G2) only. These preliminary computational results, in particular performed in a simple analytically tractable case, seem to indicate that the effect of a chronotherapy on the growth rate $\lambda(\varepsilon,\theta)$ highly depends on the amplitude $\varepsilon$ of the death rate but that the optimal phase $\theta$ (related to the best peak infusion phase) is  independent on $\varepsilon$ (see figure \ref{fig:chronotherapie1}). In future work, we intend to introduce also an effect of chronotherapy on the transition rates $K_{i\rightarrow i+1}$.

From a more general point of view (i.e., independently of chronotherapeutic considerations),
Theorem~\ref{Comparaison-locale} analytically shows, at least in the single-phase case, that under the
control of a periodic function exerting its influence on cell division, a selective advantage is given to those
cells that are able to divide with a cell cycle duration slightly lower than the control function period. But, as
numerically illustrated on Fig.~\ref{fig:Floquet-Perron_global}, this cell cycle duration should be {\em not
too much lower} than the control period, or else the advantage is lost. This leads to a biological speculation
(or prediction): in a population of proliferating cells with variable cycle duration times, all being under the
control of a common 24 h-periodic circadian clock, those cells that are well controlled by the clock, and
endowed with a cycle duration between say 21 h and 23 h should quickly outnumber the others. Hence in
proliferating healthy tissues (fast renewing tissues such as gut or bone marrow), an intrinsic cell cycle time
of 21 to 23 h should be observed (if such an observation is possible).

Now to explain the initial tumour growth data that first motivated this
study, we can speculate in the following way: tumour cells are less
sensitive than healthy cells to circadian clock control (indeed it is
known from chronotherapeutics in oncology that ``in contrast with
consistent rhythmic changes in drug tolerability mechanisms in host
tissues, tumour rhythms appear heterogeneous with regard to clock gene
expression and rhythm in pharmacology determinants as a function of
tumour type and stage"\cite{Levi}), so that their proliferation is more likely
to be governed by a simple Perron eigenvalue rather than by one of the
Floquet type. Tumour surrounding healthy cells and host immune cells, in contrast with tumour cells, are still under
circadian control and they may thus have a selective advantage over cancer cells
as long as this circadian control is present. Circadian clock disruption by
perturbed light-dark cycle destroys this advantage, and these perturbed host cells oppose in a
less efficient way local tissue invasion by cancer cells, hence the
resulting curves shown in the introduction. Of course such speculation
remains to be documented (in particular by investigating differential
circadian clock control on proliferation in tumour and healthy tissues), but
this is our best explanation so far for this phenomenon. \appendix\section{Appendix}
%
%
%
%
%
%
\subsection{Existence theory for $\lambda$}
This part is dedicated to the demonstration of the existence of the Floquet eigenvalue. Particularly, we try to prove it under general hypothesis on the periodic function $\psi$. For instance, a short adaptation of the demonstration given in \cite{MMP} would be sufficient for the case of a positive continuous periodic function $\psi$, but one would like to have the  possibility of studying non smooth functions such as a square wave (which for instance could have value $1$ during the day and $0$ during the night). We give a proof of the existence of the Floquet eigenvalue in the one-phase model. It can easily be adapted for a multiphase-model with out  death rates where the coefficients would have the form 
$K_{i\rightarrow i+1}=K_i\psi_i(t)\chi_{[a_i,+\infty[}$ with the same hypothesis on the functions $\psi_i$. We prove here existence of a solution to three eigenproblems: the direct eigenproblem
\begin{equation}
\left\lbrace\begin{array}{l}
\partial_t N(t,x)+\partial_x N(t,x)+(\lambda+K\psi(t)\chi_{[a,+\infty[}(x))N(t,x)=0,\\[0.3cm]
N(t,0)=2K\psi(t)\int_a^\infty N(t,x)dx,\\[0.3cm]
N\geq 0, \qquad \int_0^\infty Ndx=1,
\end{array}\right.
\end{equation}
the dual eigenproblem
\begin{equation}\label{adjointFloquet}
\left\lbrace\begin{array}{l}
-\partial_t \phi(t,x)-\partial_x \phi(t,x)+(\lambda+K\psi(t)\chi_{[a,+\infty[}(x))\phi(t,x)=2K\psi(t)\chi_{[a,\infty[}(x)\phi(t,0),\\[0.3cm]
\phi> 0, \quad \int_0^\infty N\phi dx=1,
\end{array}\right.
\end{equation}
and the delay differential equation
\begin{equation}
\dot P(t)=-(K\psi(t)+\lambda)P(t)+2Ke^{-\lambda a}\psi(t-a)P(t-a),\qquad P>0, \int_0^1 P(t)dt=1.
\end{equation}
We give a normalization for $P$ to ensure uniqueness.
\begin{theorem}
For any positive $T$-periodic bounded function $\psi\not=0$, $a\geq 0$
here exists a unique $\lambda,N,\phi,P$ such that $P>0$ is solution to (\ref{EDOret1}) and $N\geq 0$ is solution to (\ref{Floquet1}) ($N>0$ if $\psi$ is positive).
\end{theorem}
The proof is  based on the Krein-Rutman theorem (see \cite{Dautray-Lions} for instance).  We consider a $T$-periodic nonnegative bounded function $\psi\not=0$. We adapt the proof from \cite{MMP} to our case. First, using the methods of characteristics for the partial differential equations, we reduce the eigenproblems to integral equations on $N(t,0)$, $\phi(t,0)$ and $P$.
We consider three operators depending on a parameter $\mu$. For a bounded $T$-periodic function $\mathcal{M}$, we define $\mathcal{N}_i=\mathcal{L}_i(\mathcal{M})$ by 
\begin{eqnarray}
\mathcal{N}_1(t)&=& 2K\int_a^\infty \psi(t-x)e^{-\mu x-K\int_a^x\psi(t-x+s)ds}\mathcal{M}(t-x)dx,\\[0.3cm]
\mathcal{N}_2(t)&=& 2K\int_a^\infty \psi(t)e^{-\mu x-K\int_a^x\psi(t-x+s)ds}\mathcal{M}(t-x)dx,\\[0.3cm]
\mathcal{N}_3(t)&=& 2K\int_a^\infty \psi(t+x)e^{-\mu x-K\int_a^x\psi(t+s)ds}\mathcal{M}(t+x)dx.
\end{eqnarray}
These operators are defined such that for $\mu=\lambda$, we get,
\begin{eqnarray*}
P(t)&=& \mathcal{L}_1(P)(t),\\
N(t,0)&=& \mathcal{L}_2(N(.,0))(t),\\
\phi(t,0)&=& \mathcal{L}_3(\phi(.,0))(t).
\end{eqnarray*}
This means that the functions should be nonnegative eigenvectors of these three operators associated to the eigenvalue $1$.
\begin{lemma}
For $\mu\geq 0$, 
\begin{itemize}
	\item $\mathcal{L}_i$ maps $L^\infty_{per}(0,T)$ into itself,
	\item $\mathcal{L}_1,\mathcal{L}_3$ are continuous compact operators on $C_{per}(0,T)$,
	\item $\mathcal{L}_1,\mathcal{L}_3$ are strongly positive and $\mathcal{L}_2$ is nonnegative (strongly positive if $\psi>0$). 
\end{itemize}
\end{lemma}
\noindent {\bf Proof.}
For $\mathcal{M}$ bounded, one has, since for $x>a$, $\int_a^x \psi (t)dt\leq \langle\psi\rangle (x-a)+\|\psi\|_\infty T$, 
$$\|\mathcal{L}_i(M)\|_\infty\leq 2\frac{\|\psi\|_\infty}{\langle(\psi)\rangle}e^{K\|\psi\|_\infty T}\|\mathcal{M}\|_\infty=C\|\mathcal{M}\|_\infty.$$
For continuity and compactness we only explicit the proof for $i=1$, the case $i=3$ is very similar. We consider $\mathcal{M}$ continuous  and $h$ small,
\begin{eqnarray*}
\mathcal{N}_1(t+h)&=& 2K\int_a^\infty \psi(t+h-x)e^{-\mu x-K\int_a^x\psi(t+h-x+s)ds}\mathcal{M}(t+h-x)dx,\\
				&=& 2K\int_{a-h}^\infty \psi(t-x)e^{-\mu (x+h)-K\int_{a}^{x+h}\psi(t-x+s)ds}\mathcal{M}(t-x)dx,\\
				&=& 2K\int_{a-h}^a\psi(t-x)e^{-\mu (x+h)-K\int_{a}^{x+h}\psi(t-x+s)ds}\mathcal{M}(t-x)dx\\
				&+& 2K\int_{a}^\infty \psi(t-x)e^{-\mu (x+h)-K\int_{a}^{x+h}\psi(t-x+s)ds}\mathcal{M}(t-x)dx,\\
				&=& A_h +2K\int_{a}^\infty \psi(t-x)\bigg(e^{-\mu (x+h)-K\int_{a}^{x+h}\psi(t-x+s)ds}-e^{-\mu x-K\int_{a}^{x}\psi(t-x+s)ds}\bigg)\mathcal{M}(t-x)dx\\
				&+& 2K\int_{a}^\infty \psi(t-x)e^{-\mu x-K\int_{a}^{x}\psi(t-x+s)ds}\mathcal{M}(t-x)dx,\\
				&=& A_h +B_h + \mathcal{N}_1(t). 
\end{eqnarray*}
We have bounds on $A_h$ and $B_h$, 
\begin{equation*}\label{Ah}
|A_h| \leq  2K\|\psi\|\|\mathcal{M}\|_\infty h,
\end{equation*}
\begin{equation*}\label{Bh}
|B_h|\leq K\|\psi\|_\infty h\|\mathcal{N}_1\|_\infty\leq CK\|\psi\|_\infty \|\mathcal{M}\|_\infty h.
\end{equation*}
Therefore, using (\ref{Ah}),(\ref{Bh}), we obtain the continuity and the compactness of operator $\mathcal{L}_1$. Using the same
 techniques we can prove continuity and compactness of operator $\mathcal{L}_3$. The operator $\mathcal{L}_2$ needs regularity on $\psi$ to be
 compact (and continuous). All these operators are positive. We can apply the Krein-Rutman theorem (weak form \cite{Dautray-Lions}). We denote $\rho_1,\rho_3$ the spectral radii  of respectively $\mathcal{L}_1,\mathcal{L}_3$. They are positive (since $\mathcal{L}(1)\geq \varepsilon>0$, $\rho_1\geq\varepsilon$), so are the associated nonnegative eigenfunctions. If $\mathcal{M}_1(t)=0$, then 
 $$\psi(t-x)\mathcal{M}_1(t-x)=0,\qquad \text{for}\quad x\geq a,$$
 which leads to $\psi\mathcal{M}_1\equiv 0 $ and $\rho_1\mathcal{M}_1\equiv 0$.Therefore $\mathcal{M}_1$ and similarly $\mathcal{M}_3$ can not vanish. 

\begin{lemma}
\begin{equation*}
	\mathcal{L}_2(\psi\mathcal{M}_1)=\rho_1\psi\mathcal{M}_1, 
\end{equation*}
\begin{equation*}
	\rho_1=\rho_3.
\end{equation*}
\end{lemma} 
\noindent {\bf Proof.}
The first point is a straightforward computation. The second point uses the duality of operators $\mathcal{L}_2$ and $\mathcal{L}_3$,
\begin{eqnarray*}
\int_0^T\mathcal{L}_2(\psi\mathcal{M}_1)(t),\mathcal{M}_3(t)dt &=& \int_0^T \psi(t)\mathcal{M}_1(t)\mathcal{L}_3(\mathcal{M}_3)(t)dt,\\
\rho_1\int_0^T\psi(t)\mathcal{M}_1(t)\mathcal{M}_3(t)dt &=&\rho_3\int_0^T\psi(t)\mathcal{M}_1(t)\mathcal{M}_3(t)dt. 
\end{eqnarray*}
The existence of a solution to (\ref{EDOret1}) is equivalent to the existence of a positive fixed point of $\mathcal{L}_1$ for $\mu=\lambda$, therefore, we need to find $\mu$ such that $\rho_1(\mu)=1$. For $\mu=0$, we have 
$$\mathcal{L}_3(1)=2.$$
Therefore $\rho_1(0)=2$. As $\rho$ is a decreasing function of $\mu$ and $\rho_1(\infty)=0$, there exists some positive $\lambda$
such that $\rho_1(\lambda)=1$. The solution to (\ref{EDOret1}) is then given by such $\lambda$ and $P=\mathcal{M}_1$. Then, the function $N$ defined $N(t,0)=\psi(t)\mathcal{M}_1(t)$ and the characteristics 
$$N(t,x)=N(t-x,0)e^{-\lambda x-\int_0^x K\psi(t-x+s)\chi_{[a,\infty[}(s)ds},$$
is solution to (\ref{Floquet1}). We remark then, as $\lambda>0$, that $N(t,\infty)=0$. Similarly, we define $\phi$ by 
$$\phi(t,x)=\int_x^\infty K\psi(t+y-x)\chi_{[a,\infty[}(y)\phi(t+y-x,0)e^{-\int_x^y\lambda +K\psi(t+s-x)\chi_{[a,\infty[}(s)ds}dy.$$
This is a solution to (\ref{adjointFloquet}).

\subsection{Explicit solutions for the multiphase eigenproblem}

In the following $T=1$.

We give here explicit solutions to the eigenproblem in the multiple phases case. We do not give details for the demonstration. We consider a $3$ phase model without death terms, where the transition terms have the form:
$$K_{i\rightarrow i+1}(t,x)=K_i\psi_i(t)\chi_{[a_i,\infty[}(x).$$
Here, $\psi_i$ is a positive $1$-periodic function satisfying $\langle \psi_i\rangle =1$. We consider the following very specific case: we choose $a_1,a_2,a_3>0$ such that $a_1+a_2+a_3=1$, and we choose $\psi_i$ in the following way, for a fixed positive $1$-periodic function $\psi$,
\begin{eqnarray*}
\psi_1(t)&=&\psi(t),\\
\psi_2(t)&=&\psi(t-a_2),\\
\psi_3(t)&=&\psi(t-a_2-a_3).
\end{eqnarray*} 
To explain the form of the coefficients, we make the following remark: if we denote $P_i(t)=\int_{a_i}^\infty N_i(t,x)dx$ (the same idea as for the one phase model), the $1$-periodic functions $P_i$ satisfies a system of delay differential equations and since $a_1+a_2+a_3=1$, the $1$-periodic functions $Q_i$ defined by $Q_1(t)=P_1(t), Q_2(t)=P_2(t+a_2), Q_3(t)=P_3(t+a_2+a_3)$ satisfy a system of ordinary differential equations.
$$\frac{d}{dt}\left(\begin{array}{c}Q_1(t)\\Q_2(t)\\Q_3(t)\end{array}\right)=-\left(\begin{array}{ccc}-\lambda-K_1\psi(t)& 0 & 2K_3e^{-\lambda a_1}\psi(t)\\ K_1e^{-\lambda a_2}\psi(t) & -K_2\psi(t)-\lambda &0 \\ 0 & K_2e^{-\lambda a_3}\psi(t)& -K_3\psi(t)-\lambda\end{array}\right)\left(\begin{array}{c}Q_1(t)\\Q_2(t)\\Q_3(t)\end{array}\right).$$
We denote $M(t)$ the above matrix. Due to the special form of the functions
$\psi_i$, we have $M(t)M(t')=M(t')M(t)$, for all $t,t'$.
Therefore we can write
$$Q(t)=\exp\bigg(\int_0^tM(s)ds\bigg)Q(0).$$
The vector $Q(t)$ is $1$-periodic, thus, $Q(0)$ has to be a positive eigenvector of $\exp(\int_0^1M(s)ds)$ associated to the eigenvalue $1$.
The matrix $\exp(\int_0^1M(s)ds)$   has eigenvalue $1$ if and only if
\begin{equation}\label{lambda-analytique}
(K_1+\lambda)(K_2+\lambda)(K_3+\lambda)-2K_1K_2K_3e^{-\lambda (a_1+a_2+a_3)}=0.
\end{equation}
This leads to $Q_i(t)=e^{\lambda \int_0^t(\psi(s)-1) ds} Q_i(0)$, where $Q(0)$ is a positive vector satisfying $\int_0^1M(s)ds Q(0)=0$. Then, we can compute $P_i(t)$ and $N_i(t,0)$. Finally, using the methods of characteristics,    
the eigenfunctions $N_i$ are given, up to a normalization, by 
\begin{eqnarray*}
N_1(t,x)&=&2K_3U_3\psi(t+a_1-x)e^{\lambda\int_0^{t-x+a_1}(\psi(s)-1)ds-\lambda x-\int_{0}^x K_1\psi(t-x+s)\chi_{[a_1,\infty[}(s)ds},\\[0.3cm]
N_2(t,x)&=&K_1U_1\psi(t-x)e^{\lambda\int_0^{t-x}(\psi(s)-1)ds-\lambda x-\int_{0}^x K_2\psi(t-x+s-a_2)\chi_{[a_2,\infty[}(s)ds},\\[0.3cm]
N_3(t,x)&=&K_2U_2\psi(t-a_2-x)e^{\lambda\int_0^{t-x-a_2}(\psi(s)-1)ds-\lambda x-\int_{0}^x K_3\psi(t-x+s-a_2-a_3)\chi_{[a_3,\infty[}(s)ds},
\end{eqnarray*}
where 
$$\left(\begin{array}{c} U_1\\ U_2\\ U_3 \end{array}\right)=\left(\begin{array}{c} 1\\ \frac{K_1 e^{-\lambda a_2}}{K_2+\lambda} \\ \frac{K_1+\lambda}{2K_3e^{-\lambda a_1}}\end{array}\right).$$
The adjoint eigenfunctions are given by the formulas
\begin{eqnarray*}
\phi_1(t,x)&=&e^{-\lambda\int_0^{t-a_2-a_3-\min(x,a_1)}(\psi(s)-1)ds+\lambda \min(x,a_1)}V_1,
\\[0.3cm]
\phi_2(t,x)&=&e^{-\lambda\int_0^{t-\min(x,a_2)}(\psi(s)-1)ds+\lambda \min(x,a_2)}V_2,
\\[0.3cm]
\phi_3(t,x)&=&e^{-\lambda\int_0^{t-a_2-\min(x,a_3)}(\psi(s)-1)ds+\lambda \min(x,a_3)}V_3,
\end{eqnarray*}
where
$$\left(\begin{array}{c} V_1\\ V_2\\V_3 \end{array}\right)=\left(\begin{array}{c} 1\\ \frac{K_1+\lambda}{K_1e^{-\lambda a_1}}\\\frac{2K_3e^{-\lambda a_3}}{K_3+\lambda} \end{array}\right).$$
Basically, the ideas for the computations of $\phi_i$ are the same, based on the following remark, as 
$$\phi_i(t,x)=\int_0^\infty K_{i\rightarrow i+1}(t+y,x+y)\phi_{i+1}(t+y,0)e^{-\int_0^y \lambda+K_{i\rightarrow i+1}(t+y',x+y')dy'}dy,$$
(with a factor $2$ for $i=3$), we have $\phi_i(t,x)=\phi_(t,a_i)$ for $a\geq a_i$.  This leads to a differential equation for $\phi_i(t,0)$. Details are left to the reader.
In this case, we compute $\int_0^\infty N_i(t,x)\phi_i(t,x)dx$. As we have $\phi_i(t,x)=\phi_i(t,a_i)$ for $x\geq a_i$, 
$$\int_0^\infty N_i(t,x)\phi_i(t,x)dx=\int_0^{a_i} N_i(t,x)\phi_i(t,x)dx+\phi_i(t,a_i)\int_{a_i}^\infty N_i(t,x)dx.$$
We have 
\begin{eqnarray*}
 \int_0^{\infty}N_1(t,x)\phi_1(t,x)&=&  (K_1+\lambda)e^{\lambda a_1}\int_0^{a_1}\psi(t-x+a_1)dx+U_1V_1e^{\lambda a_1},\\[0.3cm]
 \int_0^{\infty}N_2(t,x)\phi_2(t,x)dx&=&(K_1+\lambda)e^{\lambda a_1}\int_0^{a_2}\psi(t-x)dx+U_2V_2e^{\lambda a_2},\\[0.3cm]
\int_0^{\infty}N_3(t,x)\phi_3(t,x)dx&=& (K_1+\lambda)e^{\lambda a_1}\int_0^{a_3}\psi(t-a_2-x)dx+U_3V_3e^{\lambda a_3}.
\end{eqnarray*}
Particularly, in this case, $\int_0^\infty N_i\phi dx$ is not always constant. We denote $\Psi(t)=\int_0^t(\psi(s)-1)ds$, it is a $1$ periodic function. We also denote $C_i=U_iV_ie^{\lambda a_i}, C=(K_1+\lambda)e^{\lambda a_1}$, both these constants are positive, 
\begin{eqnarray*}
 \int_0^{\infty}N_1(t,x)\phi_1(t,x)&=&  C(a_1+\Psi( t)-\Psi(t+a_1))+C_1,\\[0.3cm]
 \int_0^{\infty}N_2(t,x)\phi_2(t,x)dx&=&C(a_2+\Psi(t-a_2)-\Psi (t))+C_2,\\[0.3cm]
\int_0^{\infty}N_3(t,x)\phi_3(t,x)dx&=& C(a_3+\Psi (t+a_1)-\Psi(t-a_2))+C_3.
\end{eqnarray*}
For instance, using the parameters of the simulation, we have, $\Psi(t)=\frac{0.9}{2\pi}\sin(2\pi t)$, $a_2=0.5$,
$$\int_0^\infty N_2\phi_2(t,x)dx=(Ca_2+C_2)-2C\frac{0.9}{2\pi}\sin(2\pi t)=C'-C_2'\sin(2\pi t),$$
$$\int_0^1\gamma(t+\theta)\int_0^\infty N_2\phi_2(t,x)dx=C'\int_0^1\cos^6(2\pi(t+\theta))dt-C_2'\int_0^1\cos^6(2\pi(t+\theta))\sin(2\pi t)dt,$$
a short computation leads to 
$$\gamma(t)=\frac{1}{32}\cos(6\pi t)+\frac{3}{24}\cos(4\pi t)+\frac{15}{32}\cos(2\pi t)+\frac{5}{16},$$
therefore, 
$$\int_0^1\gamma(t+\theta)\int_0^\infty N_2\phi_2(t,x)dx=C''-C''_2\sin(2\pi\theta),$$
Where $C''\geq C_2''>0$.
Therefore, in this particular case, 
$$\lim \frac{\lambda^{\varepsilon,\theta}-\lambda^0}{\varepsilon}=C_2''\sin(2\pi\theta)-C''.$$

\textbf{Acknowledgment}
This work could not have been achieved without the help of Benoit Perthame. The authors are also grateful to Marie Doumic for very fruitful discussions.

\end{document}